\documentclass[11pt]{amsart} 

\renewcommand{\bar}{\overline}

\newcommand{\pa}{\partial}
\renewcommand{\phi}{\varphi}

\newcounter{hours}\newcounter{minutes}

\newcommand{\ka}{K\"ahler }

\font\strange=msbm10

\newcommand{\R}{{{\mathchoice  {\hbox{$\textstyle{\text{\strange R}}$}}
{\hbox{$\textstyle{\text{\strange R}}$}}
{\hbox{$\scriptstyle  N\kern-0.3em  R$}}  
{\hbox{$\scriptscriptstyle  R\kern-0.2em  R$}}}}}

\newcommand{\Z}{{{\mathchoice  {\hbox{$\textstyle{\text{\strange Z}}$}}
{\hbox{$\textstyle{\text{\strange Z}}$}}
{\hbox{$\scriptstyle  Z\kern-0.3em  Z$}}
{\hbox{$\scriptscriptstyle  Z\kern-0.2em  Z$}}}}}

\newcommand{\N}{{{\mathchoice  {\hbox{$\textstyle{\text{\strange N}}$}}
{\hbox{$\textstyle{\text{\strange N}}$}}
{\hbox{$\scriptstyle  N\kern-0.3em  N$}}
{\hbox{$\scriptscriptstyle  N\kern-0.2em  N$}}}}}

\newcommand{\bb}{{\frac{\sqrt{-1}}{2\pi}}}

\usepackage{amsmath,amsthm,amscd}    
\usepackage{calc} 

\title 
[]{A note on the holomorphic invariants of Tian-Zhu}
\author{Zhiqin Lu}
\date{May 24, 2001}

\address[Zhiqin Lu]
{Department of Mathematics\\ 
University of California, Irvine\\ 
Irvine, CA 92697} 
\email[Zhiqin Lu]{zlu@math.uci.edu}
\thanks{Research supported by NSF grant DMS 0196086}

\newtheorem{theorem}{Theorem}

\newtheorem{lemma}{Lemma}
\newtheorem{cor}{Corollary}

\theoremstyle{remark}

\begin{document}
\maketitle




In this short note, we compute the holomorphic invariants
defined by Tian and Zhu~\cite{TZ} on smooth hypersurfaces of
$CP^n$. The holomorphic invariants, which  generalize 
the famous Futaki invariants~\cite{Fu1},
are obstructions
towards the existence of K\"ahler-Ricci solitons.

For a K\"ahler manifold with the first positive Chern class, the
existence of the K\"ahler-Ricci soliton can be reduced to the
existence of the solution of a non-linear equation of Monge-Ampere
type. In general, solving such an equation is highly non-trivial. 
Similar to the Futaki invariants, the Tian-Zhu invariants gives the
obstruction {\sl before} one need to solve the equation. 
It is thus very important to compute it concretary. In this paper, 
in the case of hypersurfaces, we give an explicit formula.

Let $M\subset CP^n$ be a smooth hypersurface defined by a
homogeneous polynomial 
$F=0$ of degree $d$. Let $v$ and $X$ be two 
holomorphic vector fields
on
$CP^n$. For the sake of simplicity, we assume that
\[
v=\sum_{i=0}^n v^iZ_i\frac{\pa}{\pa Z_i},\,\,\, {\rm and}\quad
X=\sum_{i=0}^n X^iZ_i\frac{\pa}{\pa Z_i},
\]
where $[Z_0,\cdots,Z_n]$  is the homogeneous coordinate of
$CP^{n}$, $(v^0,\cdots,v^n)\in C^{n+1}$, $(X^0,\cdots,X^n)\in
C^{n+1}$. We  further assume
that
\begin{equation}
\sum_{i=0}^n v^i=0, \quad\sum_{i=0}^n X^i=0.
\end{equation}
If $v$ and $X$ are tangent vector fields of $M$, then
there are complex numbers $\lambda$ and $\kappa$ such that
\begin{equation}\label{1}
vF=\kappa F,\quad XF=\lambda F. 
\end{equation}
Let $\omega$ be the \ka form of the Fubini-Study metric of $CP^n$.
Then $(n-d+1)\omega$ restricts to a representative of the first
Chern class $c_1(M)$ of $M$.  
Thus there is a smooth function $\xi$ on $M$ such that
\[
{\rm Ric}((n-d+1)\omega|_M)-(n-d+1)\omega|_M=\pa\bar\pa\xi.
\]
For fixed holomorphic vectors $X$ and $v$,
the holomorphic invariant defined by Tian-Zhu~\cite{TZ}, in our
context, is
\begin{equation}\label{2}
F_X(v)=(n-d+1)^{n-1}\int_M v(\xi-(n-d+1)\theta_X)e^{(n-d+1)\theta_X
}\omega^{n-1},
\end{equation}
where $\theta_X$ is defined as
\begin{equation}\label{5}
\left\{
\begin{array}{l}
i(X)\omega=\bb\bar\pa\theta_X,
\\
\int_M e^{(n-d+1)\theta_X}\omega^{n-1}=d.
\end{array}
\right.
\end{equation}

The main property of the Tian-Zhu invariants is the following
(cf. ~\cite{TZ}):

\begin{theorem}
Let $F_X(v)$ be the Tian-Zhu invariant. Then we have

1. If the K\"ahler-Ricci soliton exists, that is, we have
\[
{\rm Ric}(\omega)-\omega=L_X\omega
\]
for some K\"ahler metric $\omega$. Then 
$F_X(v)\equiv 0$.

2. $F_X(v)$ is independent of the choice of the K\"ahler metric
$\omega$ within  the first Chern class.
\end{theorem}

In this note, we give a ``computable'' expression of $F_X(v)$. Our
main result is as follows: 
\begin{theorem}\label{main}
Using the notations as above,
defined the function
\begin{equation}\label{10}
\phi(X)=\sum_{k=0}^\infty\frac{n!(n-d+1)^k}{(n+k)!}
\sum_{\alpha_0+\cdots+\alpha_n=k}
X_0^{\alpha_0}\cdots X_n^{\alpha_n},
\end{equation}
where $\alpha_0,\cdots,\alpha_n\in \Z^{n+1}$ are nonnegative
integers. Let
\begin{equation}\label{11}
\sigma(X)=(-\frac{\lambda(n-d+1)}{n}+d)\phi(X)+\frac{d}{n}
\sum_{i=0}^n X^i\frac{\pa\phi(X)}{\pa X^i}.
\end{equation}
Then the invariants defined by Tian-Zhu can be explicitly
expressed as
\begin{equation}\label{12}
F_X(v)=-(n-d+1)^{n-1}d\left(\kappa+\sum_{i=0}^n
v^i\frac{\pa\log\sigma(X)}{\pa X^i}\right).
\end{equation}
\end{theorem}

\begin{cor}
The Futaki invariant for the hypersurface $M$ is
\[
F(v)=-(n-d+1)^{n-1}\frac{(n+1)(d-1)}{n}\kappa.
\]
\end{cor}
\qed

The rest of this note is devoted to the proof Theorem~\ref{main}.
We define
\begin{equation}
\tilde\theta_X=\frac{\lambda_0|Z_0|^2+\cdots+\lambda_n|Z_n|^2}{
|Z_0|^2+\cdots+|Z_n|^2}.
\end{equation}
Then we have
\begin{equation}
i(X)\omega=\bar\pa\tilde\theta_X.
\end{equation}
By comparing the above equation with~\eqref{5}, we have
\begin{equation}\label{21}
\theta_X=\tilde\theta_X+c_X
\end{equation}
for a constant $c_X$.
First, we have the following lemma

\begin{lemma}\label{p}
\[
\int_{CP^n}e^{(n-d+1)\tilde\theta_X}\omega^n=\phi(X),
\]
where $\phi(X)$ is defined in~\eqref{10}.
\end{lemma}

{\bf Proof.} This follows from the expansion 
\[
e^{(n-d+1)\tilde\theta_X}=\sum_{k=0}^\infty
\frac{(n-d+1)^k}{k!}\tilde\theta_X^k,
\]
and the elementary Calculus.

\qed

\begin{lemma}\label{q}
Using the same notation as above, we have
\[
F_X(v)=(n-d+1)^{n-1}\left(-\kappa
d-\int_M(n-d+1)\theta_ve^{(n-d+1)\theta_X}
\omega^{n-1}\right).
\]
\end{lemma}

{\bf Proof.}
By~\cite[Theorem 4.1]{Lu9}, we have
\[
{\rm div}\,v+v(\xi)+(n-d+1)\theta_v=-\kappa,
\]
where $\theta_v$ is the function on $CP^n$ defined by
\[
\theta_v=\frac{v_0|Z_0|^2+\cdots+v_n|Z_n|^2}
{|Z_0|^2+\cdots+|Z_n|^2},
\]
and $\kappa$ is defined in~\eqref{1}.
Then
~\eqref{2} becomes 
\begin{align}\label{4}
\begin{split}
&\qquad F_X(v)=(n-d+1)^{n-1}\\
&\cdot\left(\int_M(-\kappa-{\rm div}\,
v-(n-d+1)\theta_v
-(n-d+1)v(\theta_X))e^{(n-d+1)\theta_X}\omega^{n-1}\right).
\end{split}
\end{align}
We also have
\begin{equation}\label{3}  
{\rm div}\,(e^{(n-d+1)\theta_X}v)
=
e^{(n-d+1)\theta_X}({\rm div}\, v+(n-d+1)v(\theta_X)).
\end{equation}
The lemma follows from ~\eqref{5}, ~\eqref{4}, ~\eqref{3}
and the divergence theorem.

\qed

The following key lemma transfers the integration on $M$ to the
integrations on $CP^n$.

\begin{lemma}\label{r}
\begin{equation}\label{13}
(n-d+1)\int_{M}\theta_v e^{(n-d+1)\theta_X}\omega^{n-1}
=d\sum_{i=0}^n v^i\frac{\pa\log \sigma}{\pa X^i},
\end{equation}
where $\sigma(X)$ is defined in~\eqref{11}.
\end{lemma}

{\bf Proof.}
Let
\begin{equation}\label{8} 
\eta=\log\frac{|F|^2}{(|Z_0|^2+\cdots+|Z_n|^2)^d}.
\end{equation}
Then $\eta$ is a smooth function on $CP^n$ outside $M$. 
We have the following identity:
\begin{align}
\begin{split}
&
\bar\pa(e^{(n-d+1)\theta_X}\pa\eta\wedge\omega^{n-1})-\frac
{n-d+1}{n}
i(X)(e^{(n-d+1)\theta_X}\pa \eta\wedge\omega^{n})\\&
=-e^{(n-d+1)\theta_X}\pa\bar\pa\eta\wedge\omega^{n-1}-\frac
{n-d+1}{n}
e^{(n-d+1)\theta_X}(\lambda-d\tilde\theta_X)\omega^n.
\end{split}
\end{align}
Since on $CP^n$, there are no $(2n+1)$ forms, the left hand side 
of the above equation is the divergence of some vector field. 
Integrate the equation  on both side and use the divergence theorem,
we have
\begin{equation}\label{14}
\int_{CP^n}e^{(n-d+1)\theta_X}\pa\bar\pa\eta\wedge\omega^{n-1}
=-\frac{n-d+1}{n}
\int_{CP^n}(\lambda-d\tilde\theta_X)e^{(n-d+1)\theta_X}\omega^n.
\end{equation}
     
By~\cite[page 388]
{GH}, in the sense of currents, we have
\begin{equation}\label{9}
\pa\bar\pa\eta=[M]-d\omega.
\end{equation}
Thus from ~\eqref{14},
\begin{align}\label{15}
\begin{split}&
\int_M e^{(n-d+1)\theta_X}\omega^{n-1}
=(-\frac{\lambda(n-d+1)}{n}+d)\int_{CP^n}
e^{(n-d+1)\theta_X}\omega^n\\&\qquad+\frac{d(n-d+1)}
{n}\int_{CP^n}\tilde\theta_X
e^{(n-d+1)\theta_X}\omega^n.
\end{split}
\end{align}
From Lemma~\ref{p}, we have
\begin{equation}\label{17}
\sum_{i=0}^nX^i\frac{\pa\phi(X)}{\pa X^i}
=(n-d+1)\int_{CP^n}\tilde\theta_Xe^{(n-d+1)
\tilde\theta_X}\omega^n.
\end{equation}
By~\eqref{21}, ~\eqref{15}and~\eqref{17}
\begin{equation}\label{18}
\int_M e^{(n-d+1)\theta_X}\omega^{n-1}=\sigma(X)
e^{c_X}.
\end{equation} 
From the above equation, we have
\begin{equation}\label{19}
(n-d+1)\int_M \theta_ve^{(n-d+1)\theta_X}\omega^{n-1}
=\sum_{i=0}^n v^i\frac{\pa\sigma(X)}{\pa X^i}e^{c_X}.
\end{equation}
On the other hand, from~\eqref{18}, we have
\begin{equation}\label{20}
d=\sigma(X) e^{c_X},
\end{equation}
by~\eqref{5}.
Lemma~\ref{r} follows from~\eqref{19} and~\eqref{20}.

\qed

Theorem~\ref{main} follows from Lemma~\ref{q} and Lemma~\ref{r}.

\bibliographystyle{abbrv}
\bibliography{bib}

\end{document}